\documentstyle[12pt,righttag,ctagsplt]{amsart}
\newtheorem{prop}{Proposition}
\newtheorem{th}[prop]{Theorem}
\newtheorem{cor}[prop]{Corollary}
\newtheorem{lem}[prop]{Lemma}

\theoremstyle{definition}
 
\theoremstyle{definition}
\newtheorem{example}{Example}
\theoremstyle{remark}
 
\newcommand{\bbN}{{\Bbb{N}}}
\newcommand{\bbR}{{\Bbb{R}}}
\newcommand{\bbZ}{{\Bbb{Z}}}

\newcommand{\al}{\alpha}

\newcommand{\de}{\delta}
\newcommand{\e}{{\varepsilon}}

\newcommand{\f}{\varphi}

\renewcommand{\span}{\operatorname{span}}
\newcommand{\supp}{\operatorname{supp}}
\newcommand{\rank}{\operatorname{rank}}
\newcommand{\card}{\operatorname{card}}
\newcommand{\codim}{\operatorname{codim}}

\newcommand{\sgn}{\operatorname{sgn}}

\newcommand{\lb}{\label}

\newcommand{\lra}{\longrightarrow}

\newcommand{\wtw}{if and only if }
\newcommand{\buo}{without loss of generality }

\newcommand{\ONTO}{\buildrel {\mbox{\small onto}}\over \longrightarrow} 
\newcommand{\DEF}{\buildrel {\mbox{\rm\small def}}\over =}

\makeatletter
\def\@currentlabel{2.1}\label{e:dispaa}
\def\@currentlabel{2.21}\label{e:dispau}
\def\@currentlabel{2.22}\label{e:dispav}
 \def\@currentlabel{2.23}\label{e:dispaw}
\def\@currentlabel{2.24}\label{e:dispax}
\makeatother

\makeatletter
\def\alphenumi{%
  \def\theenumi{\alph{enumi}}%
  \def\p@enumi{\theenumi}%
  \def\labelenumi{(\@alph\c@enumi)}}
\makeatother
\font\tenex=cmex10
\newskip\ttglue

\def\eightpoint{\def\rm{\fam0\eightrm}
  \textfont0=\eightrm \scriptfont0=\sixrm \scriptscriptfont0=\fiverm
  \textfont1=\eighti  \scriptfont1=\sixi  \scriptscriptfont1=\fivei
  \textfont2=\eightsy  \scriptfont2=\sixsy  \scriptscriptfont2=\fivesy
\textfont3=\tenex  \scriptfont3=\tenex  \scriptscriptfont3=\tenex
\textfont\itfam=\eightit  \def\it{\fam\itfam\eightit}
\textfont\slfam=\eightsl  \def\sl{\fam\slfam\eightsl}
\textfont\ttfam=\eighttt  \def\tt{\fam\ttfam\eighttt}
\textfont\bffam=\eightbf  \scriptfont\bffam=\sixbf
\scriptscriptfont\bffam=\fivebf  \def\bf{\fam\bffam\eightbf}
\tt  \ttglue=.5em plus.25em minus.15em
\normalbaselineskip=9pt
\setbox\strutbox=\hbox{\vrule height7pt depth2pt width0pt}
\let\sc=\sixrm  \let\big=\eightbig \normalbaselines\rm}
\font\eightrm=cmr8 \font\sixrm=cmr6 \font\fiverm=cmr5
\font\eighti=cmmi8  \font\sixi=cmmi6   \font\fivei=cmmi5
\font\eightsy=cmsy8  \font\sixsy=cmsy6 \font\fivesy=cmsy5
\font\eightit=cmti8  \font\eightsl=cmsl8  \font\eighttt=cmtt8
\font\eightbf=cmbx8  \font\sixbf=cmbx6 \font\fivebf=cmbx5

\def\eightbig#1{{\hbox{$\textfont0=\ninerm\textfont2=\ninesy
	\left#1\vbox to6.5pt{}\right.\enspace$}}}

\begin{document}
\title{One-complemented subspaces of real sequence spaces }
\author{Beata Randrianantoanina}
\address{Department of Mathematics and Statistics \\ Miami University \\
	 Oxford, OH 45056}
\email{randrib@@muohio.edu}


\begin{abstract}
Characterizations are given for  1-complemented hyperplanes of 
strictly monotone real
Lorentz spaces and 1-complemented finite codimensional 
subspaces (which
contain at least one basis element) of real Orlicz 
spaces equipped with
either Luxemburg or Orlicz norm.
\end{abstract}
\maketitle

\section{Introduction}

The question about the form of projections and the structure of their
ranges  arises naturally  in geometry
of Banach spaces and has many applications in other fields. Thus it has
been studied by many authors since the time of Banach,
 but still even the structure of
norm-one projections and 1-complemented subspaces is far from understood
for most Banach spaces (for more detailed
discussions of applications and some of the open questions 
see, e.g., the surveys \cite{ChP,Doust}).

In this paper we study which subspaces of finite codimension are
1-complemented in real Orlicz and Lorentz sequence spaces.  In the earlier
paper \cite{complex}, the author gave necessary conditions for
1-complementability of subspaces of complex spaces with 1-unconditional
bases.  However, the methods of the complex case do not transfer to the real
case.  The main difference is that every 1-complemented subspace of a
complex space with 1-unconditional basis also has a 1-unconditional basis
(\cite{KW,F83}), but this fact fails for subspaces of real spaces
(\cite{L79,BFL}), and even for subspaces of symmetric real spaces (since
every space with 1-unconditional basis is 1-complemented in some symmetric
space \cite{L72}, \cite[Theorem~3.b.1]{LT1}).

Among the real sequence spaces, only 1-complemented subspaces of $\ell_p$,
\ $1 \leq p \leq \infty$ are fully understood.  Namely, it is well known
that a subspace $Y \subset \ell_p$, $ (1 \leq p < \infty, \ p \neq 2)$ is
1-complemented if and only if $Y$ is spanned by a block basis of a
permutation of the original basis (\cite[Theorem~2.a.4]{LT1}).  The only
other available description is a sufficient condition for
1-complementability of subspaces of symmetric sequence spaces; namely
subspaces spanned by a block basis with constant coefficients of a
permutation of the original basis are 1-complemented in symmetric spaces
(\cite[Proposition~3.a.4]{LT1}).

In the present paper, we show that the  sufficient condition above is also
necessary for 1-complementability of: 
\begin{enumerate}
\item hyperplanes in strictly
monotone Lorentz spaces not isometric to $\ell_p, \ 1 \leq p < \infty$,
(Corollary~\ref{corLor}), and 
\item finite codimensional subspaces containing
at least one basis element of Orlicz spaces (with either Luxemburg or Orlicz
norm) not isomorphic to any $\ell_p, \ 1 \leq p < \infty$
(Corollary~\ref{orlicz}).
\end{enumerate}

We use a technique of numerically positive operators (defined by
formula~\eqref{numpos} below), and our main tool is that a projection $P$
has norm one if and only if $(Id-P)$ is numerically positive
(Proposition~\ref{A}).

We use standard Banach space notation; all undefined terms may be
found, for example, in \cite{LT1}.

\section{Preliminaries}

Throughout the paper we will consider a real Banach space $X$ with
1-unconditional basis $\{e_i\}_{i\in I}$. Most often $X$ will be either a
Lorentz or an Orlicz sequence space.  We recall their definitions below
(see \cite{LT1}).

Let $1 \leq p < \infty$ and let $w = (w_n)_n$ be a non-increasing sequence
of nonnegative numbers with $w_1 = 1$.  Then the {\it Lorentz space}
$\ell_{w, p}$ is a space of all sequences $x = (x_n)_{n}$ for which $\| x
\|_{w, p} < \infty$ where
\begin{equation*}
\| x \|_{w, p} = \left( \sum^{\dim \ell_{w, p}}_{n=1} w_n \left( x^*_n
\right)^p \right)^{\frac 1p}
\end{equation*}
and where $\left( x^*_n \right)_n$ is the non-increasing rearrangement of
$( | x_n |)_n$.

A function $\varphi : [0, \infty) \longrightarrow [0, \infty)$ is called an
{\it Orlicz function} if $\varphi$ is non-decreasing, convex and
$\varphi(1)= 1$.  An {\it Orlicz space} $\ell_\varphi$ is a space of
sequences $x = (x_n)^{\dim \ell_\varphi}_{n=1}$ such that there exists
$\lambda > 0$ with

\begin{equation*}
\sum^{\dim \ell_\varphi}_{n=1} \varphi \left( \frac{| x_n |}\lambda \right) <
\infty \ .
\end{equation*}
It is customary to consider two (equivalent) norms on $\ell_\varphi$ : the
Luxemburg and the Orlicz norm.  The Luxemburg norm is given by

\begin{equation*}
\| x \|_\varphi = \inf \left\{ \lambda > 0 : \sum^{\dim \ell_\varphi}_{n=1}
\varphi \left( \frac{| x_n |}\lambda \right) \leq 1 \right\}
\end{equation*}
and the Orlicz norm is given by
\begin{equation*}
\| x \|_{\varphi, O} = \sup \left\{ \sum^{\dim \ell_\varphi}_{n=1} x_n y_n
: \sum^{\dim \ell_\varphi}_{n=1} \varphi^* (| y_n |) \leq 1 \right\} \ ,
\end{equation*}
where $\varphi^*$ denotes the Orlicz function complementary to $\varphi$ in
the sense of Young, i.e. for any $u \geq 0$
\begin{equation*}
\varphi^* (u) = 
\sup \{ tu - \varphi (t) : 0 < t < \infty \} \ .
\end{equation*}

Let $X$ be a real Banach space.
We say that a functional $x^* \in
X^*$ is {\it norming} for $x \in X$ \wtw $\| x^* \| = \| x \|$ and
$x^*(x) = \| x \|^2$.
Denote $\Pi(X) = \{ (x,x^*) : \|x\| = 1 \text{ and
} x^*
\text{ is norming for x } \} \subset (X,X^*)$.
Following Rosenthal
\cite{R84} we say that an operator
$T:X\lra X$ is
{\it numerically positive}  if
\begin{equation}\lb{numpos}
x^*(Tx)\ge 0
\end{equation}
 whenever
$(x,x^*)\in\Pi(X).$
This is equivalent to requiring a slightly weaker
condition that given $x$ with $\|x\|=1$ there exists
$x^*$ so that
$(x,x^*)\in\Pi(X)$ and $x^*(Tx)\ge 0$
(see Lumer \cite{L61,L63}).
  By results of Lumer \cite{L61} and Lumer and Phillips
\cite{LP}  (see also \cite{BD}) it is equivalent to
the requirement
that
$\|\exp(-\alpha T)\|\le
1$ for $\alpha\ge 0.$
This immediately implies the following simple fact ,
which is the main tool for our applications :
\begin{prop}  $(\text{see} $ \cite{R84,KR}$)$ \label{A}
If $X$ is a real Banach space and  \/ $P:X\lra X$ is a
projection
 then $P$ is numerically positive if and only if \/
$\|Id - P\|=1,$ where $Id$ denotes the identity on $X$.
\end{prop}

We need to introduce the following technical property.
We say that a space $X$ with 1-unconditional basis
$\{e_i\}_{i\in I}$
has property $(P)$ if for all $i,j\in I$ and for all $\e >0$
$$\|e_i +\e e_j\|_X > \|e_i\|_X =1\ .$$

Notice that property $(P)$ is slightly weaker than strict
monotonicity, e.g.,
a Lorentz sequence spaces $\ell_{w,p}$ satisfies $(P)$ whenever $w_2\ne0$,
where  $w=(w_1,w_2,w_3,\ldots)$ for all $p,1\le p<\infty$, i.e.,
$\ell_{w,p} \in (P)$ whenever $\ell_{w,p}$ is not isometric to $\ell_\infty$.

An Orlicz sequence space $\ell_\varphi$ satisfies $(P)$ whenever
$\varphi (t)>0$ for all $t>0$ and $\varphi (1)=1$.

We say that a space $X$ with 1-unconditional basis $\{e_i\}_{i\in I}$ has
property $(Q)$ if for all $i,j\in I$
$$\lim_{\e\to0} {\|e_i+\e e_j\|_X-1 \over \e} = 0\ .$$

Notice that a Lorentz sequence space $\ell_{w,p}$ satisfies $(Q)$ whenever
$p>1$, and Orlicz sequence space $\ell_\varphi$ satisfies $(Q)$ whenever
$\varphi'(0)=0$.

 We will be interested in finite codimensional subspaces of real
Banach spaces.
Notice that a subspace $Y\subset X$,
$\codim Y = n$, is spanned by disjointly supported vectors \wtw
$Y$ can be presented in the form 
$Y= \bigcap_{j=1}^n \ker g_j$, where functionals $g_j$ are such
that $\card(\supp g_j)\le2$ for all $j\le n$.
We will consider the ``standard'' form of $Y\subset X$,
$\codim Y = n$, i.e.,
$Y= \bigcap_{j=1}^n \ker f_j$, where functionals $f_j$ are such
that $f_{jk} = \delta_{jk}$ for all $j,k\le n$.
Further, we will consider a projection
$P:X \ONTO Y$,
where  $P = Id-\sum_{j=1}^n f_j\otimes u_j$;
here $Id$ denotes the identity on $X$ and $u_1,\ldots,u_n$
are linearly independent elements
in $X$ with $f_j(u_k)=\delta_{jk}$.
In this notation we have:

\begin{lem}\label{supp}
Suppose that $X$ satisfies properties $(P)$ and $(Q)$ and that the
projection
$P:X\to Y$ has norm one.
Then
$$\bigcup_{j=1}^n \supp u_j = \bigcup_{j=1}^n \supp f_j\ .$$
\end{lem}

\begin{pf}
If $k\notin \bigcup_{j=1}^n \supp u_j$ then
$$(P(e_k))_k = 1-\sum_{j=1}^n f_j (e_k) u_{jk} =1$$
Since $\|P\|=1$, property $(P)$ implies that $P(e_k)= e_k$.
Hence $e_k\in \bigcap_{j=1}^n \ker f_j=Y$ and
$k\notin \bigcup_{j=1}^n \supp f_j$.
Therefore $\bigcup_{j=1}^n \supp f_j \subset \bigcup_{j=1}^n \supp u_j$.

For the other inclusion, consider $k\notin \bigcup_{j=1}^n \supp f_j$.
Let $M= \max\{ |u_{jj}| : j=1,2,\ldots,n\}+1$.
Fix arbitrary $i$, $1\le i\le n$, $\e \ne0$, and consider element
$x_{\e} = e_k + \e e_i$.
Denote by $x_{\e}^N = a_{\e} e_k^* + b_{\e} e_i^*$  a norming element
for $x_{\e}$ with $\|x_{\e}^N\|_{X^*} =1$.
By property $(P)$ of $X$ $b_{\e}\ne0$ whenever $\e \ne0$.
Since, by Proposition~\ref{A},
 $\sum_{j=1}^n f_j\otimes u_j$ is numerically positive, we get
for all $\e \ne0$:
\begin{eqnarray*}
0&\le & \sum_{j=1}^n f_j (x_{\e} ) x_\e^N (u_j)
= \sum_{j=1}^n (f_{jk} + \e f_{ji}) (a_\e u_{jk} + b_\e u_{ji}) \\
& = & \sum_{j=1}^n \e\delta_{ji} (a_\e u_{jk} + b_\e u_{ji})
= \e (a_\e u_{ik} + b_\e u_{ii}) \ .
\end{eqnarray*}
Thus for all $\e \ne0$
\begin{equation}\lb{supp1}
\sgn (a_\e u_{ik} + b_\e u_{ii})  = \sgn\e = \sgn (b_\e)\ .
\end{equation}
By property $(Q)$ of $X$ $\lim_{\e \to0} b_\e =0$.
Thus \eqref{supp1} implies that $u_{ik} =0$.
Indeed, if $|u_{ik}| = \eta >0$ let $\e$ be small enough so that
$|b_\e| <{\eta\over 2M}$  and $a_\e >{1\over2}$.
Then
$$|b_\e u_{ii}| < {\eta\over 2M} \cdot M
= {\eta \over 2} < a_\e \cdot |u_{ik}|\ .$$
Thus $\sgn (a_\e u_{ik} + b_\e u_{ii}) = \sgn (a_\e u_{ik})$ for all
$\e$, contradicting \eqref{supp1}.
\end{pf}

\section{Lorentz sequence spaces $\ell_{w,p}$, $p>1$}

In this section we study 1-complemented hyperplanes of $\ell_{w,p}$
where $p>1$.

\begin{th} \lb{Lor}
Let $\ell_{w, p}$ be a Lorentz sequence space with $1 < p < \infty$ and
$w_2 > 0$.  Suppose that $Y = \ker f$ is 1-complemented in $\ell_{w, p}$
and
$\card(\supp f) \geq n > 2$.  Then $p=2$ and $1=w_1=w_2\ldots=w_n$.
\end{th}

Our formulation resembles the characterizations of 1-complemented
finite codimensional subspaces of $\ell_p$ in terms of their representation
as an intersection of kernels of functionals $Y=\bigcap \ker f_j$ (see
\cite{BP88}).

We postpone the proof of Theorem~\ref{Lor} to the end of this section.

As a consequence of Theorem~\ref{Lor}, we obtain the following
characterization of $\ell_2$ among Lorentz sequence spaces:

\begin{cor} \lb{lorentzl2}
Let $\ell_{w, p}$ be a Lorentz sequence space with
$1 < p < \infty$ and
$w_i > 0$ for all $i\ge1$.
 Suppose that there exists a functional $f$ such that
$ \card(\supp f) > 2$ and  $\ker f$ is 1-complemented in
$\ell_{w, p}$.  Then
$\ell_{w, p} = \ell_2$.
\end{cor}

In the above corollary, the assumption about nonzero weights is equivalent
to asking that $\ell_{w, p}$ be strictly monotone. It can
be slightly relaxed if $\card (\supp f) = \infty$ or if 
$\dim \ell_{w, p} <
\infty$.\begin{itemize}
\item[(a)] If $\card (\supp f) = \infty$ it is enough to assume that 
$w_2 > 0$.
\item[(b)] If $\dim \ell_{w, p} = d < \infty$ and 
$\card (\supp f) = n \leq d$ it is
enough to assume that $w_{d - n+3} > 0$.
\end{itemize}
However the assumption about nonzero weights cannot be completely removed
as the following example demonstrates.
\begin{example}
Let $w=(1,1,1,0)$ and consider a 4-dimensional space $\ell^4_{w,2}$.
Let $Y=\ker(e_1^* +e_2^*+e_3^*)$ and $P_2:\ell_2^3\lra\ell_2^3$ be the
orthogonal projection. Define $P:\ell^4_{w,2}\lra\ell^4_{w,2}$ by
$$P(x_1,x_2,x_3,x_4) = P_2(x_1,x_2,x_3) + 0\cdot e_4.$$

Then:
$$\|P(x_1,x_2,x_3,x_4)\|_{w,2} = \|P_2(x_1,x_2,x_3)\|_2 
\le \|(x_1,x_2,x_3,0)\|_2 \le\|(x_1,x_2,x_3,x_4)\|_{w,2} .$$

Thus $\|P\|\le1$ i.e. $\ker(e_1^* +e_2^*+e_3^*)$ is 1-complemented
in $\ell^4_{w,2}\ne \ell^4_{2}$.
\end{example}

The condition $\card (\supp f) \le 2$ means that $\ker f$ is 
spanned by disjointly supported
vectors. In \cite{complex} the author studied the form of 1-complemented
disjointly spanned subspaces of Lorentz spaces and she obtained:

\begin{th}(\cite[Theorem~6.3]{complex})
\label{glorentz}
Let $\ell_{w,p}$, with $1 < p < \infty$, be a real or complex Lorentz
sequence space.  Suppose that $\{ x_i \}_{i \in I}$ are mutually
disjoint elements of $\ell_{w,p}$ such that $\card (I) \geq 2$ and $F
= \overline{\span} \{ x_i \}_{i \in I} $ is $1$-complemented in
$\ell_{w,p}$.  Suppose, moreover, that $w_\nu \neq 0$ for all $\nu
\leq \Sigma \DEF \sum_{i \in I} \card (\supp x_i)\quad (\leq \infty)$.

Then
\begin{itemize}
\item[$(a)$] $w_\nu = 1$ for all $ \nu \leq \Sigma$,
\end{itemize}

\noindent
or
\begin{itemize}
\item[$(b)$] $ |x_{il}| = |x_{ik}|$ for
all $i \in I$ and all $k, l \in \supp x_i$.
\end{itemize}
\end{th}

Thus, as an immediate consequence of Theorem~\ref{Lor} and
\cite[Theorem~6.3]{complex}, we get:

\begin{cor}\lb{corLor}
Let $\ell_{w, p}$ be a Lorentz sequence space with $1 < p < \infty$ and
$w_k > 0$ for all $k$, i.e. $\ell_{w, p}$ is strictly monotone.  
Suppose that $Y = \ker f$ is 1-complemented in
$\ell_{w, p}$ and $\card (\supp f) = 2$, i.e. $f = f_i {e_i}^* + f_j
{e_j}^*$ for some $i \neq j$.

Then $|f_i| = |f_j|$ for all $i, j$, 
or $\ell_{w, p} = \ell_p$, i.e. $w_k = 1$ for all $k$.
\end{cor}

\begin{pf*}{Proof of Corollary~\ref{lorentzl2}}

If $2 < \card (\supp f) = \dim \ell_{w, p} \leq \infty$ then the corollary
follows immediately from Theorem~\ref{Lor}.

Assume thus, that $2 < \card (\supp f) = n < \dim \ell_{w, p} \leq \infty$.
We will prove by induction that $w_m = w_1 = 1$ for all $m \leq \dim
\ell_{w, p}$.  The first induction step follows from Theorem~\ref{Lor}, also we
get that $p = 2$.

Assume now that $w_m = 1$ for some $n \leq m < \dim \ell_{w, p}$.  We will
show that $w_{m+1} = 1$.

Let $\{ i_1, \ldots, i_{m - n+1} \} \subset \{ 1, \ldots, \dim \ell_{w, 2}
\}$ be such that $\{ i_1, \ldots, i_{m - n+1} \} \cap \supp f =
\emptyset$.

Consider a Lorentz space $\ell_{w' , 2} = \span\{ e_i : i \notin \{
i_1, \ldots, i_{m-n+1} \}\}$, where $w'$ is a weight defined by
$w'_k = w_{k+m-n+1}$.

By inductive hypothesis, $w'_1 = \ldots = w'_{n-1} = w_m = 1$.

For any $x \in \ell_{w, p}$ with $\supp x \subseteq \supp f$, consider an
element
\begin{equation*}
x_+ \DEF x + \sum^{m-n+1}_{k=1} \| x \|_{w, 2} e_{i_k} \ .
\end{equation*}
Since $\{ e_{i_1}, \ldots, e_{i_{m-n+1}} \} \subset \ker f$, we have
\begin{equation*}
Px_+ = Px + \sum^{m-n+1}_{k=1} \| x \|_{w, 2} Pe_{i_k} = Px +
\sum^{m-n+1}_{k=1} \| x \|_{w, 2} e_{i_k} \ .
\end{equation*}
Since $P$ is a contractive projection $\| Px \|_{w, 2} \leq \| x \|_{w, 2}$,
and thus $\| x \|_{w,2} \geq |(Px)_i|$ for all
$i \in \supp x$. Clearly also $\|x\|_{w,2} \geq |x_i|$ for all
$i \in \supp x$.
Hence
\begin{equation*}
\|x_+\|^2_{w,2} = \sum^{m-n+1}_{k=1} \|x\|^2_{w, 2} + \|x\|^2_{w',2}
\end{equation*}
and
\begin{equation*}
\| Px_+ \|^2_{w, 2} = \sum^{m-n+1}_{k=1} \|x\|^2_{w, 2} + \|Px\|^2_{w',2} \ .
\end{equation*}
Since $P$ is contractive in $\ell_{w, 2}$ we conclude that $\|Px\|_{w
', 2} \leq \|x\|_{w',2}$ for all $x$ with $\supp x \cap \{ i_1,
\ldots, i_{m-n+1} \} = \emptyset$.  That is, $P$ is contractive in 
$\ell_{w', 2}$.  Thus $\ker f$ is 1-complemented in $\ell_{w',2}$ and,
by Theorem~\ref{Lor}, since $\card (\supp f) = n$ and $w'_2 = w_{m-n+3}
\ge w_m > 0$, we get $w_n ' = 1$.
Thus $w_{m+1} = 1$, as claimed.
\end{pf*}

\begin{pf*}{Proof of Theorem~\ref{Lor}}
Since $\ell_{w,p}$ is symmetric we can assume
without loss of generality that
$$f= (f_1,f_2,f_3,\ldots,f_n,\ldots)$$
where $f_1\ge f_2\ge\dots\ge f_n >0$.
Then $f\otimes u$ is numerically positive
$(u= (u_1,u_2,\ldots,u_n,\ldots))$
and, by Lemma~\ref{supp}, $u_1,u_2,\ldots,u_n\ne0$.

We organize the proof of Theorem~\ref{Lor} into four assertions:

\noindent
{\bf{Assertion 1:}}
If $n \geq 4$ or $n=3$ and $f_3 < f_1$ then $p=2$ and $$\frac{w_2 + \ldots +
w_{n-1}}{n-2} \sum_{i=1}^{n-2} u_i = \frac{f_1 + \ldots + f_{n-2}}{f_n} u_n$$
$$u_n = \frac{f_n}{f_{n-1}} w_n u_{n-1}$$

\noindent
{\bf{Assertion 2:}}
If $n \geq 4$ or ($n=3$ and $f_2 < f_1$) then $p=2$ and $w_n = 1$, i.e.
Theorem~\ref{Lor} holds.

\noindent
{\bf{Assertion 3:}}
If $n=3$ and $f_1 = f_2$ then $p=2$, $w_2 = 1$ and
$$u_1 = \frac{f_1}{f_3}
w_3 u_3.$$

\noindent
{\bf{Assertion 4:}}
If $n=3$ then $p=2$ and $w_3 = 1$, i.e. Theorem~\ref{Lor} holds.

In the proof of all assertions, we will use the following elements $x(a)$
and $x(a, \e)$:
\begin{equation*}
x(a) = e_1 + \ldots + e_{n-2} - \frac{f_1 + \ldots + f_{n-2}}{f_{n-1}} a
e_{n-1} 
- \frac{f_1 + \ldots + f_{n-2}}{f_n} (1-a) e_n 
\end{equation*}
\begin{align*}
x(a, \e) \ = &x(a) - \e e_{n-1} - \e e_n \cr
= \ &e_1 + \ldots + e_{n-2} - \biggl( 
\frac{f_1 + \ldots + f_{n-2}}{f_{n-1}} a
+ \e \biggr) e_{n-1} - \biggl( \frac{f_1 + \ldots + f_{n-2}}{f_n} (1-a) + 
\e \biggr) e_n \ ,
\end{align*}
where $a \in [0, 1], \ \e \in [-1, 1]$.

Notice that for all $a , \ x(a) \in \ker f$.

For the proof of Assertion 1, let $$\eta = \min \biggl( \frac{f_{n-1}}{f_1
+ \ldots + f_{n-2}} \ , 1 - \frac{f_n}{f_1 + \ldots + f_{n-2}} \biggr).$$

Notice that by assumptions of Assertion 1, $\eta > 0$.

For any $a$ with $0 < a < \eta$ we have
\begin{equation}\lb{npf1} 
\frac{f_1 + \ldots + f_{n-2}}{f_n} (1-a) > 1 > \frac{f_1 + \ldots +
f_{n-2}}{f_{n-1}} > 0 \ .
\end{equation}
For any $a$ with $0 < a < \eta$ define
$$\delta (a) = \min \biggl( \frac{f_1 + \ldots + f_{n-2}}{f_n} (1-a) - 1, \
1 - \frac{f_1 + \ldots + f_{n-2}}{f_{n-1}}, \ \frac{f_1 + \ldots +
f_{n-2}}{f_{n-1}} \biggr) \ .$$

Then, for any $\e$ with $| \e | < \delta (a)$ 
$$| x (a, \e)_n | > | x
(a, \e)_1 | = \ldots = | x (a, \e)_{n-2} | > | x (a, \e)_{n-1} | \ .$$

Thus, a norming functional for $x (a, \e)$ can be chosen as follows:
\begin{align}
x(a, \e)^N = & \sum_{i=1}^{n-2} \biggl( \frac{w_2 + \ldots + w_{n-1}}{n-2}
\biggr) {e_i}^* - \biggl( \frac{f_1 + \ldots + f_{n-2}}{f_{n-1}} a + \e
\biggr)^{p-1} w_n e_{n-1}^* \cr
& - \biggl( \frac{f_1 + \ldots + f_{n-2}}{f_n} (1-a) + \e \biggr)^{p-1}
{e_n}^* \lb{norming} \ .
\end{align}
Since $f \otimes u$ is numerically positive, for all $a, \ 0 < a < \eta$
and $\e$ with $| \e | < \delta (a)$ we have
$$f \bigl( x (a, \e) \bigr) x(a, \e)^N (u) \geq 0 \ .$$

Hence:
\begin{align*}
f \bigl( x (a, \e) \bigr) = \ &f \bigl( x(a) \bigr) - \e f (e_{n-1} + e_n)
\cr
= \ &0 - \e (f_{n-1} + f_n) = - \e (f_{n-1} + f_n) \ . 
\end{align*}
\begin{align*}
x (a, \e)^N (u) = \ & \frac{w_2 + \ldots + w_{n-1}}{n-2} \sum_{i=1}^{n-2}
u_i - \biggl( \frac{f_1 + \ldots + f_{n-2}}{f_{n-1}} a + \e \biggr)^{p-1} w_n
u_{n-1} \cr
&\qquad - \biggl( \frac{f_1 + \ldots + f_{n-2}}{f_n} (1-a) + \e
\biggr)^{p-1} u_n \ .
\end{align*}

Thus $$- \sgn \e = \sgn \biggl( f \bigl( x (a, \e) \bigr) \biggr) = \sgn
\biggl( x (a, \e)^N (u) \biggr) \ .$$
Since $p > 1$ we conclude that 
\begin{equation}
x (a, \e)^N  (u) = 0 \lb{step1} \ .
\end{equation}
Hence:
\begin{align}\lb{gl}
\frac{w_2 + \ldots + w_{n-1}}{n-2} \sum_{i=1}^{n-2} u_i -& \biggl(
\frac{f_1 + \ldots + f_{n-2}}{f_{n-1}} \biggr)^{p-1} w_n u_{n-1} a^{p-1}
\cr
&\qquad - \biggl( \frac{f_1 + \ldots + f_{n-2}}{f_n} \biggr)^{p-1} u_n
(1-a)^{p-1} = 0 \ . 
\end{align}

Notice that \eqref{gl} is valid for all $a$ with $0 < a < \eta$ and all
parameters in \eqref{gl} except $a$ are fixed, i.e. \eqref{gl} in fact
is:
$$A + Ba^{p-1} + C (1-a)^{p-1} = 0$$ where $A, \ B, \ C$ are constants such
that $A, \ C \neq 0$; and $B = 0$ if and only if $w_n = 0$.

If $B = 0$ we get $$(1-a)^{p-1} = \biggl( - \frac AC \biggr) = const.$$
for all $a \in (0, \eta)$, which is impossible since $p > 1$.

If $B \neq 0$ we differentiate and get
$$\biggl( \frac{a}{1-a} \biggr)^{p-2} = - \frac CB = const.$$
Thus $p = 2$.

When $p = 2$, equation \eqref{gl} becomes
\begin{align*}
& \biggl( \frac{w_2 + \ldots + w_{n-1}}{n-2} \sum_{i=1}^{n-2} u_i -
\frac{f_1 + \ldots + f_{n-2}}{f_n} u_n \biggr) \cr
& \qquad + a \biggl( \frac{f_1 + \ldots + f_{n-2}}{f_n} u_n - \frac{f_1 +
\ldots + f_{n-2}}{f_{n-1}} w_n u_{n-1} \biggr) = 0
\end{align*}
for all $a \in (0, \eta)$.

Thus 
\begin{equation}
\frac{w_2 + \ldots + w_{n-1}}{n-2} \sum_{i=1}^{n-2} u_i = \frac{f_1
+ \ldots + f_{n-2}}{f_n} u_n \lb{gl1}
\end{equation}
and 
\begin{equation}
u_n = \frac{f_n}{f_{n-1}} w_n u_{n-1} \ . \lb{gl2}
\end{equation}
This finishes the proof of Assertion 1.

For the proof of Assertion 2, assume that $n \geq 4$ or 
that ($n=3$ and $f_2
< f_1$).

Notice first that the assumptions of Assertion 2 
are stronger than those of
Assertion 1.  Thus, we have $p=2$.

Consider elements $x (a, \e)$ with $a=1$ and 
$$| \e | < {\e_1} \DEF \min
\biggl( \frac{f_1 + \ldots + f_{n-2}}{f_{n-1}} - 1, \ 1 \biggr) \ .$$
Notice that by our assumptions, $\e_1 > 0$.

Then 
$$x(1, \e) = e_1 + \ldots + e_{n-2} - 
\biggl( \frac{f_1 + \ldots +
f_{n-2}}{f_{n-1}} + \e \biggr) e_{n-1} - \e e_n$$

and $$| x(1, \e)_{n-1} | > | x(1, \e)_1 | = \ldots = | x(1, \e)_{n-2} | > |
x(1, \e)_n | \ .$$

Thus a norming functional for $x(1, \e)$ can be chosen to be
$$x(1, \e)^N = \biggl( \frac{w_2 + \ldots + w_{n-1}}{n-2} \biggr)
\sum_{i=1}^{n-2} {e_i}^* - \biggl( \frac{f_1 + \ldots + f_{n-2}}{f_{n-1}} +
\e  \biggr) e_{n-1}^* - \e w_n {e_n}^* \ .$$

Similarly as in \eqref{step1} we conclude that $$x(1, 0)^N (u) = 0 \ .$$

Hence $$\frac{w_2 + \ldots + w_{n-1}}{n-2} \sum_{i=1}^{n-2} u_i = \frac{f_1
+ \ldots + f_{n-2}}{f_{n-1}} u_{n-1} \ .$$

Combining this with \eqref{gl1} we get $$u_n = u_{n-1}$$ and further, by
\eqref{gl2}, $$\frac{f_n}{f_{n-1}} w_n = 1 \ .$$

But $f_n \leq f_{n-1}$, thus $$w_n = \frac{f_{n-1}}{f_n} \geq 1 \ .$$

Hence $w_n = 1$ and Assertion 2 is proved.

For the proof of Assertion 3, assume that $n=3$ and $f_1 = f_2$.
\newline
By symmetry of $\ell_{w, p}$ we can assume without loss of generality that
$u_1 = u_2$.
Similarly as in previous cases, consider the element 
$x(a, \e)$ with $a
\in ( {f_2}/({f_2 + f_3}),  1 )$ and 
$$| \e | < \e(a)
\DEF \min \biggl( 1-a, \ \frac 12  \bigl( a - \frac{f_1}{f_3} (1-a)
\bigr) \biggr) \ .$$
Notice that $\e (a) > 0$ when $a \in ( {f_2}/({f_2 + f_3}), 1)$. Then
$$x(a, \e) = e_1 - (a + \e) e_2 - \biggl( \frac{f_1}{f_3} (1-a) + \e
\biggr) e_3$$
and $$| x(a, \e)_1 | > | x(a, \e)_2 | > | x(a, \e)_3 | \ .$$

Thus the norming functional for $x(a, \e)$ is given by
$$x(a, \e)^N = {e_1}^* - (a + \e)^{p-1} w_2 {e_2}^* - \biggl(
\frac{f_1}{f_3} (1-a) + \e \biggr)^{p-1} w_3 {e_3}^* \ .$$
Similarly as in \eqref{step1}, we conclude that $$x(a, 0)^N (u) = 0 \ .$$
Thus, for all $a \in ( {f_2}/({f_2 + f_3}), 1 ),$
\begin{equation}
u_1 - a^{p-1} w_2 u_2 - \biggl( \frac{f_1}{f_3} (1-a) \biggr)^{p-1} w_3
u_3 = 0 \ . \lb{gl3}
\end{equation}

By the same argument as in \eqref{gl}, we see that $p=2$; and \eqref{gl3}
becomes
$$u_1 - \frac{f_1}{f_3} w_3 u_3 + a \biggl( \frac{f_1}{f_3} w_3 u_3 - w_2
u_2 \biggr) = 0 \ .$$

Since this equality holds for all $a \in ( {f_2}/({f_2 + f_3}), 1)$
we have 
\begin{equation}
u_1 = \frac{f_1}{f_3} w_3 u_3 \lb{as3}
\end{equation}

and $$\frac{f_1}{f_3} w_3 u_3 = w_2 u_2 \ .$$
Since $u_1 = u_2$ we conclude that $w_2 = 1$, which ends the 
proof of Assertion~3.

Finally, for the proof of Assertion 4, assume that $n=3$.

If $f_2 < f_1$ then by Assertion 2, Theorem~\ref{Lor} holds.

If $f_2 = f_1$ and $f_3 < f_1$ then $u_1 = u_2$ and by Assertion 3, $p=2, \
w_2 = 1$ and $$u_1 = \frac{f_1}{f_3} w_3 u_3 \ .$$
On the other hand, by Assertion 1, $$u_3 = \frac{f_3}{f_2} w_3 u_2 \ .$$
Thus 
$$u_1 = \frac{f_1}{f_3} w_3 \frac{f_3}{f_2} w_3 u_2 = {w_3}^2 u_1 \ .$$
Hence $w_3 = 1$ and Theorem~\ref{Lor} holds.

If $f_1 = f_2 = f_3$ then by symmetry of $\ell_{w, p}$, \ $u_1 = u_2 = u_3$
and by Assertion 3, $$u_1 = \frac{f_1}{f_3} w_3 u_3 = w_3 u_3 = w_3 u_1 \
.$$

Thus $w_3 = 1$ and Theorem~\ref{Lor} holds.
\end{pf*}

\section{One-complemented subspaces of Orlicz sequence spaces}

We say that $\f$  is  
{\it similar} to $t^p$ for some $p\in[1,\infty)$ if there exist $C, t_0>0$
so that $\f(t) = Ct^p$ for all $t<t_0$. 
We say that $\f$   is {\it equivalent} to $t^p$ for some $p\in[1,\infty)$
if there exist $C_1, C_2, t_0>0$
so that $C_1t^p\le\f(t) \le C_2t^p$ for all $t<t_0$.
 
The main result of this section is the following: 

\begin{th}\lb{or} 
Let $\ell_\varphi$ be an Orlicz space such that $\f$  is not 
similar to $t^2$ and $\f(t)>0$ for all $t>0.$ Consider $\ell_\f$
with either Luxemburg or Orlicz norm.
Let $F\subset \ell_\varphi$ be a subspace of 
codimension~$n$ with $\dim F>1$.  
If $F$ contains at least one basis vector and 
$F$ is 1-complemented in $\ell_\varphi$ then $F$ can be 
represented as $F= \bigcap_{j=1}^n \ker f_j$ where $\card (\supp f_j)\le2$ 
for all $j$.

Moreover, if $\varphi$ is not equivalent to $t^p$ for any 
$p\in[1,\infty)$,
then 
$|f_{ji}|$ is either $1$ or $0$ for all $i,j$. 

If $\varphi$ is equivalent but not similar 
to $t^p$ for some $p\in[1,\infty)$,  then there 
exists 
$\gamma\ge1$ such that 
$$\{|f_{ji}| : j\le n\ ,\ i\le \dim \ell_\varphi\} \subset 
\{\gamma^m : m\in \bbZ \} \cup \{0\}\ .$$ 
\end{th}

Similarly as in the case of Lorentz spaces
our formulation resembles the characterizations of 1-complemented
finite codimensional subspaces of $\ell_p$ in terms of their representation
as an intersection of kernels of functionals $Y=\bigcap \ker f_j$ (see
\cite{BP88}).

As an immediate consequence of Theorem~\ref{or} 
we obtain the following:

\begin{cor}\lb{coror}
Let $\ell_\varphi$ be an Orlicz space with either Luxemburg or Orlicz norm,
such that $\f$  is not 
equivalent to $t^p$ 
for any $1\le p\le\infty$. 
Let $F\subset \ell_\varphi$ be a subspace of codimension~$n$ with
$\dim F>1$ 
and containing at least one basis vector. 
Then $F$ is 1-complemented in $\ell_\varphi$ 
if and only if $F$ is a span 
of a block basis with constant coefficients.
\end{cor} 

Before the proof of Theorem~\ref{or} we note that the assumptions that
$\f$ is is not 
similar to $t^2$ and $\f(t)>0$ for all $t>0$ cannot be removed.
Indeed in \cite{complex} we presented an example of an Orlicz function
$\f$ with $\f(t) = t^2$ for all $t\le a$, where $\sqrt{2/3}<a<1$, and such 
that $\ker(e_1^* +e_2^*+e_3^*)$ is 1-complemented in $\ell_\f$ 
(\cite[Example~3]{complex}). This example can be easily modified so that
for arbitrarily small $a>0$ we have $\f(t) = t^2$ when $t\le a$, and
$\ell_\f$ contains a 1-complemented subspace of codimension $n$ which
cannot be represented as a span of disjointly supported vectors.

\begin{example}
Let $\f$ be an Orlicz function such that $\f(t) = t^2$ for all $t$ 
with $0\le t\le a$ for some $a>0$. Let $m\in\bbN$ be such that $1/m\le a^2$.
Consider $\ell_\f^{3m}$ and a subspace $F\subset\ell_\f^{3m}$ defined by:
$$ F=\{(x_i)_{i=1}^{3m}\ :
\ x_1 +x_2 +x_3 = 0, x_k = x_{k+3j}, k=1,2,3, j=1,\dots,m-1\}.$$
Then $F$ cannot be presented in the form described in the conclusion
of Theorem~\ref{or}, but  $F$ is 1-complemented in $\ell_\f^{3m}$.

Note first that if $x\in F$ then $\|x\|_\f =\|x\|_2$. 
Indeed:
$$1 = \sum_{i=1}^{3m} \f\left(\frac{|x_i|}{\|x\|_\f}\right) = 
m \sum_{k=1}^{3} \f\left(\frac{|x_k|}{\|x\|_\f}\right).$$

Thus for all $i\le3m$ we have:
$$ \f\left(\frac{|x_i|}{\|x\|_\f}\right)\le\frac1m\le a^2.$$
Since $\f$ is increasing 
$$\f\left(\frac{|x_i|}{\|x\|_\f}\right) = \frac{|x_i|^2}{\|x\|_\f^2}$$
for all $i\le3m$ i.e. $\|x\|_\f =\|x\|_2$.

Let $Q:\ell_2^{3m}\lra F$ be a contractive projection of $\ell_2^{3m}$ onto
 $F$. Then $Q$ is also contractive when considered as a projection from 
$\ell_\f^{3m}$ onto
 $F$. Indeed for all $x\in\ell_\f^{3m}$ we have:
$$\|x\|_\f \ge\|x\|_2\ge\|Qx\|_2=\|Qx\|_\f.$$

\end{example}

We were unable to eliminate the condition that subspace $F$ contains 
at least one basis vector in the assumptions of Theorem~\ref{or}.
However, it follows from the author's earlier work that this 
condition is satisfied in Orlicz spaces which are either $p$-convex
with constant 1 for some $2< p < \infty$, or $q$-concave
with constant 1 for some $1< q < 2$ and $\f$ is smooth at 1, provided
that $\dim F <1/2\dim X$. Namely, we have:

\begin{th} (\cite[Theorem~1]{Rseq}) \lb{pro}
Let  $X$ be a
strictly monotone
 sequence space $(\dim X = d \ge 3)$ with a 1-unconditional basis
$\{e_i \}_{i=1}^d.$ Suppose that

$(a)$ \  $X$  is  $p$-convex with constant 1, $2<p<\infty,$

\noindent or

$(b)$ \  $X$ is
$q$-concave with constant 1, $1<q<2,$  and smooth at each basic vector.

Then any
1-complemented subspace $F$ of codimension $n$ in $X$ contains all
 but at most $2n$ basic
vectors of $X.$
\end{th}

Recall (see \cite[Definition 1.d.3]{LT2})  that a Banach space $X$ is 
{\it $p$-convex with constant 1} (resp.
{\it  $q$-concave  with constant 1}) 
 if
for every
choice of elements $\{ x_i \}_{i=1}^n $ in $X$ the following inequality
holds:
\begin{align*}
   \| \left( \sum_{i=1}^n |x_i|^p \right)^{1/p} \| \ \ &\le \ \ \left(
\sum_{i=1}^n
\|x_i\|^p
\right)^{1/p} \ \ \ \ \ \ \ \hfill \text{if } 1\le p < \infty, \\
\intertext{or, respectively,}
   \| \left( \sum_{i=1}^n |x_i|^q \right)^{1/q} \| \ \ &\ge \ \ \left(
\sum_{i=1}^n
\|x_i\|^q
\right)^{1/q} \ \ \ \ \ \ \ \hfill \text{if } 1\le q < \infty.
\end{align*}

\begin{pf*}{Proof of Theorem~\ref{or}} 
We start with choosing a convenient notation. 
Let $F= \bigcap_{j=1}^n \ker f_j$ where $f_j$'s are in the 
reduced form, 
i.e., $f_{ji} = \delta_{ij}$ for all $i,j\le n$ and $f_{j,n+1}=0$ for 
all $j$. 
If $\dim \ell_\varphi \le n+2$ there is nothing to prove so we will assume 
that $\dim \ell_\varphi \ge n+3$. 

We claim that if there exist $j\le n$ and $i,k>n$ with $f_{ji},f_{jk}\ne0$ 
then $\f$  is  
similar to $t^2$ or $\f(t)\not>0$ for all $t>0.$

Assume, for contradiction, that $F$ is 1-complemented and, say, 
$f_{1,n+2},
f_{1,n+3}\ne0$. 
Then $P : \ell_\varphi \ONTO F$ has form 
$$P(x) = x- \sum_{j=1}^n f_j (x) u_j,$$ 
where vectors $\{u_j\}_{j=1}^n$ are linearly independent. 
It is clear that we can select a subset $S'$ of indices such that 
$\{u_j|_{S'}\}_{j=1}^n$ are linearly independent. 
Let $S= \{1,2,\ldots,n+3\}\cup S'$. 
For convenience of notation we can assume \buo that $S= \{1,\ldots,M\}$, 
where $M\le \min \{\dim\ell_\varphi ,2n+3\}$. 

Let 
$$L= \Big\| \sum_{j=1}^n \biggl( \sum_{i=n+2}^M |f_{ji}| +1\biggr) 
e_j + \sum_{j=n+2}^M e_j \Big\|_\varphi\ .$$ 
Here we denote by $\|\cdot\|_\f$ either the
Luxemburg
or Orlicz norm in $\ell_\f$.

For any $k\in \{1,\ldots,n\}$, $\e \in\bbR$ with  $|\e|<1$ and any sequence
$(\alpha_j)_{j=n+2}^M \subset [-1,1]$, consider  elements 
$y$ and $y{(k,\e)}$ in $X$ 
defined by 
$$y = \sum_{j=1}^n \biggl( -\sum_{i=n+2}^M f_{ji}\alpha_i\biggr) 
e_j  +\sum_{j=n+2}^M \alpha_j  e_j\in\bigcap_{j=1}^n \ker f_j ,$$ 
$$y{(k,\e)} = y+ \e e_k.$$
 Then $\|y{(k,\e)}\|_\varphi \le L$. 
Thus there exists $\alpha_{k,\e}$ so that 
$\|y{(k,\e)} + \alpha_{(k,\e)} e_{n+1}\|_\varphi =L$. 
Set $x{(k,\e)} = y{(k,\e)} + \alpha_{k,\e} e_{n+1}$ and 
$x=y(1,0) +\alpha_{1,0} e_{n+1}$. 
Then $x\in \bigcap_{j=1}^n \ker f_j$ and 
$f_j (x{(k,\e)})=\de_{jk}\e$  for all  
$1\le k, j\le n$.

Since $Id-P = \sum_{j=1}^n f_j \otimes u_j$ is numerically positive, we 
conclude that 
\begin{equation}\lb{or0}
\sgn (x{(k,\e)}^N (u_k)) = \sgn \e
\end{equation}
where $x{(k,\e)}^N$ denotes a norming functional for $x{(k,\e)}$. 
By \cite{GrzH} and \cite{ChHK} in both cases of Luxemburg
and Orlicz norm on $\ell_\f$, a norming functional $x({k,\e})^N$ may be
 defined by the following formula:
$${x(k,\e)^N}_j = {1\over C} \sgn (x(k,\e)_j)\varphi' 
\left( {|x{(k,\e)_j}| \over L}\right) \ ;\qquad j=1,\ldots,n$$
Here $\varphi'$ denotes the left derivative of $\varphi$ and $C>0$ is a 
constant depending on $x({k,\e})$, (recall that $L= \|x_{(k,\e)}\|_\varphi$ 
for all $(k,\e)$; of course $C$ depends also on the choice of either the
Luxemburg
or Orlicz norm). 

Notice that for almost all $(\alpha_j)_{j=n+2}^M \in [-1,1]^{M-n-1}$, 
$\varphi'$ is continuous at all points 
$\{ {|x_{(k,0),j}|\over L}\ : \ {j=1}\dots M\}$. 
Denote the set of all such $(\alpha_j)_{j=n+2}^M$ by $\Lambda$. When 
$(\alpha_j)_{j=n+2}^M\in\Lambda$, since for every $k$,$\ \ $  $x(k,0) =x$, 
equation \eqref{or0} implies that for 
all $k$ with 
$1\le k\le n$
$$x^N (u_k)=0\ ,$$ 
that is 
\begin{equation}\lb{or1}
\sum_{j=1}^n - \sgn \left( \sum_{i=n+2}^M f_{ji} \alpha_i\right) 
\varphi' \left( {|\sum\limits_{i=n+2}^M f_{ji}\alpha_i| \over L}\right) 
u_{kj} + 
\sum_{j=n+2}^M \sgn (\alpha_j) \varphi' 
\left( {|\alpha_j|\over L}\right) u_{kj} =0\ . 
\end{equation}

Thus, for any $(M-n-1)$-tuple $(\alpha_j)_{j=n+2}^M$ in $\Lambda$, 
\eqref{or1} can be treated as an equation of variables $\{u_{kj}\}$ 
and it has 
$n$ linearly independent solutions: 
$$u_1|_S, u_2|_S,\ldots, u_n|_S\ .$$

Now we will look at equation \eqref{or1} with specially selected
$(M-n-1)$-tuples $(\alpha_j)_{j=n+2}^M$.
 For $r=1,\ldots, M-n-1$ consider $(M-n-1)$-tuples
$(\al_j^r)_{j=n+2}^M\in\Lambda$ such that 
\begin{equation*}
\alpha_j^r =
\begin{cases} 0 \ \ \ \ \ {\text when } \ \ \ j\ne n+r,\\ 
\delta \ \ \ \ \ {\rm when } j= n+r,
\end{cases}
\end{equation*}
where $1/2<\delta \le1$ is chosen so that 
all tuples $(\alpha^r_j)_{j=n+2}^M$ are in $\Lambda$. 
Further let 
$$\alpha_j^{M-n} = \left\{ 
\begin{array}{ll} 
0&\mbox{if } j>n+3\\ 
a&\mbox{if } j=n+2\\
b&\mbox{if } j=n+3
\end{array}
\right.  $$ 
where $a,b\in [-1,1]^2$ are such that 
$(\alpha_j^{M-n+1})_{j=n+2}^M\in\Lambda$. 
Denote by $\Lambda'$ the set of all admissable pairs $(a,b)$.
Notice that $\Lambda'\subset [-1,1]^2$ has full measure. 

Now we can consider a system of $(M-n)$ equations with $(M-1)$ variables, each 
of the form \eqref{or1} with coefficients determined by the above
defined $(M-n-1)$-tuples 
$(\alpha_j^r)_{j=n+2}^M$ for $r=1,\ldots,M-n$. 
That is, we consider a system of equations with the
matrix $A$ whose $r$-th row is given by
\begin{equation*}
\begin{split}
[A_r] &= [-\sgn (f_{1,n+r})\varphi'\left( 
{\delta|f_{1,n+r}|\over L}\right),\ldots, -\sgn(f_{n,n+r})
\varphi'\left( {\delta |f_{n,n+r}|\over L}\right), \\
      &\qquad\underbrace{0,\dots,0}_{(r-1)},
\varphi' \left( {\delta\over L}\right),
0,\ldots,0],
\end{split}
\end{equation*}
when $1\le r \le M-n-1$, and 
\begin{equation*}
\begin{split}
 [A_{M-n}] &= [-\sgn (af_{1,n+2} +bf_{1,n+3})\varphi'
\left( {|af_{1,n+2} +bf_{1,n+3}|\over L}\right),\ldots,\\
&\quad\qquad -\sgn (af_{n,n+2} +bf_{n,n+3})\varphi' 
\left( {|af_{n,n+2} +bf_{n,n+3}|\over L}\right),
\sgn(a)\varphi'\left({|a|\over L}\right),\\
&\quad\qquad \sgn(b)\varphi'\left({|b|\over L}\right),
0,\ldots,0].
\end{split}
\end{equation*}

Clearly $\rank A \ge M-n-1$. 
Also, since $u_1|_S,\ldots,u_n|_S$ are linearly independent solutions, 
the solution space of $A$ has dimension greater or equal than $n$. 
Thus $\rank A\le M-n-1$. 

Hence $\rank A= M-n-1$. Since rows $[A_r]$ for $1\le r\le M-n-1$
are clearly linearly independent, hence
 the last row of $A$ is the linear combination
of the first $(M-n-1)$ rows. 
Therefore, for all $(a,b)\in\Lambda'$ we have 
\begin{equation}\lb{or2} 
\begin{split}
 &{\sgn (a) \varphi'\left( {|a|\over L}\right) \over 
	\varphi' \left( {\delta\over L}\right) } 
	\left(-\sgn (f_{1,n+2})\varphi' 
	\left( {\delta |f_{1,n+2}|\over L}\right) \right) +\\
 &\qquad + {\sgn (b) \varphi'\left( {|b|\over L}\right) \over
	\varphi' \left( {\delta\over L}\right) } 
	\left(-\sgn (f_{1,n+3}) \varphi' 
	\left( {\delta |f_{1,n+3}|\over L}\right) \right) =\\
&\qquad = -\sgn (af_{1,n+2} +bf_{1,n+3})\varphi' 
	\left( {|af_{1,n+2} + bf_{1,n+3}| \over L}\right)
\end{split}
\end{equation}
In particular, when $b$ approaches $0$ (and respectively when
 $a$ approaches $0$) we obtain that $\lim_{t\to0}\f'(t) = 0$ and
\begin{eqnarray}
&&{-\sgn (af_{1,n+2})\varphi'\left({|a|\over L}\right)\over
	\varphi' \left( {\delta\over L}\right) } 
	\varphi' \left( {\delta |f_{1,n+2}|\over L}\right) 
	= -\sgn (af_{1,n+2})\varphi' 
	\left( {|af_{1,n+2}|\over L}\right) 
	\lb{or2'}\\
&&{-\sgn (bf_{1,n+3})\varphi'\left({|b|\over L}\right)\over 
	\varphi' \left( {\delta\over L}\right) } 
	\varphi' \left( {\delta |f_{1,n+3}|\over L}\right) 
	= -\sgn (bf_{1,n+3})\varphi' 
	\left( {|bf_{1,n+3}|\over L}\right) 
	\lb{or2''}
\end{eqnarray}

Therefore, if $\sgn (a)=\sgn (f_{1,n+2})$ and 
$\sgn (b) = \sgn (f_{1,n+3})$ and if we denote 
$$x = {af_{1,n+2} \over L}\quad ,\quad 
y = {bf_{1,n+3}\over L}$$ 
then \eqref{or2'}, \eqref{or2''} and \eqref{or2} imply that: 
\begin{equation}\lb{or3} 
\varphi' (x+y) = \varphi' (x) + \varphi'(y)
\end{equation}
and \eqref{or3} is valid for almost all $(x,y)\in [0,{1\over L}] \times 
[0,{1\over L}]$. 
Since $\varphi'$ is increasing we conclude that $\varphi'$ is linear on 
$[0,{2\over L}]$ and thus there exists $C\ge0$ so that 
$$\varphi (t)=Ct^2 \ \ \ \ \ {\rm for }\ \ \ \ 
t\le {2\over L}.$$ 

But then, if $C=0$ then there exists $t>0$ with $\varphi(t)=0$, or, if
$C>0$ then $\f$ is similar to $t^2$,  which contradicts our assumptions. 
Hence $f_{1,n+2},f_{1,n+3}$ cannot both be nonzero.

The last part of Theorem~\ref{or} follows from the following result:

\begin{th} (\cite[Theorem~6.1]{complex})
\label{orlicz}
Let $\ell_\phi$ be a (real or complex) Orlicz space and let $x,y \in
\ell_\phi$, be disjoint elements such that $\|x\|_\phi = \|y\|_\phi =
1$ and $\span \{x,y \} $ is $1$-complemented in $\ell_\phi$.

Then one of three possibilities holds:
\begin{itemize}
\item[(1)] $\card (\supp x) < \infty$ and $|x_i| = |x_j|$
for all $i,j \in \supp x$; or
\item[(2)] there exists $p$, $1 \leq p \leq \infty$, such that
$\phi(t) = Ct^p$ for all $t \leq \|x\|_\infty$; or
\item[(3)] there exists $p$, $1 \leq p \leq \infty$, and constants
$C_1, C_2, \gamma \geq 0$ such that
$C_2 t^p \leq\nobreak \phi(t) \leq C_1 t^p$
for all $t \leq \|x\|_\infty$
and  such that, for all $j \in \supp x$,
$$|x_j| = \gamma^{k(j)} \cdot \|x\|_\infty$$ for some
$k(j) \in \bbZ$.
\end{itemize}
\end{th}

\end{pf*}

  \bibliographystyle{plain}
  \bibliography{tref}

\end{document}